\documentclass [11pt] {amsart}
\usepackage{amsmath,amsthm,amsfonts,epsf,amssymb}
\usepackage{epsfig, psfrag, subfigure, bm}
\newtheorem {theorem}{Theorem} [section]
\newtheorem {definition}{Definition} 
\newtheorem {corollary}{Corollary} 
\newtheorem {lemma}{Lemma}[section]
\newtheorem {result}{Result}[section] \newtheorem {observation}{Observation}
\newtheorem {conjecture}{Conjecture} 
\newtheorem {remark}{Remark} \newtheorem {example}{Example} 
 
\newcommand{\bthrm}{\begin{theorem}} \newcommand{\ethrm}{\end{theorem}}
\newcommand{\bdf  }{\begin{definition}} \newcommand{\edf}{\end{definition}}
\newcommand{\bcor }{\begin{corollary}} \newcommand{\ecor}{\end{corollary}}
\newcommand{\bres }{\begin{result}} \newcommand{\eres}{\end{result}}
\newcommand{\blem }{\begin{lemma}} \newcommand{\elem}{\end{lemma}}
\newcommand{\bconj}{\begin{conjecture}} \newcommand{\econj}{\end{conjecture}}
\newcommand{\brem }{\begin{remark}} \newcommand{\erem}{\end{remark}}
\newcommand{\bex  }{\begin{example}} \newcommand{\eex}{\end{example}}
\newcommand{\beq  }{\begin{equation}} \newcommand{\eeq}{\end{equation}}
\newcommand{\bea  }{\begin{eqnarray}} \newcommand{\eea}{\end{eqnarray}}
\newcommand{\ben  }{\begin{enumerate}} \newcommand{\een}{\end{enumerate}}
\newcommand{\bpf  }{\begin{proof}} \newcommand{\epf}{\end{proof}}
\newcommand{\bobs}{\begin{observation}} \newcommand{\eobs}{\end{observation}}
\newcommand{\be}[2]{\begin{#1} #2 \end{#1}}

\newcommand{\m}{\mathbb }

\def\HF{\widehat{HF}}
\def\HFK{\widehat{HFK}}
\def\CFK{\widehat{CFK}}

\def\->{\rightarrow}
\def\vv{\nu}
\def\ee{\varepsilon}
\def\rr{\rho}
\numberwithin{equation}{section}

\newcommand{\calh}{\mathcal{H}} 
\begin{document}
\title [A modification of the Sarkar-Wang algorithm]
{A modification of the Sarkar-Wang algorithm and an analysis of its computational complexity}
\author [Jonathan  Hales]{Jonathan  Hales}
\author [Dmytro Karabash]{Dmytro Karabash}
\author [Michael T. Lock]{Michael T. Lock}
\maketitle

\begin{abstract}
The Sarkar-Wang algorithm computes the hat version of the Heegaard Floer homology of a
closed oriented three manifold.  This paper analyzes the computational complexity of the Sarkar-Wang algorithm; then the algorithm is modified to obtain a better bound.  Then the computational complexity of calculating $\HFK$ from a Heegaard diagram by means of the modified Sarkar-Wang algorithm is also analyzed.  Under certain assumptions it is shown that the modified Sarkar-Wang algorithm is faster than the Manolescu-Ozsv\'{a}th-Sarkar algorithm  \end{abstract}

\tableofcontents

\section{Introduction}

The Heegaard Floer homology is an invariant of a closed 3-manifolds that was developed by Peter Ozsv\'{a}th and Zolt\'{a}n Szab\'{o} \cite{Introduction}. The knot version of this invariant, known as the Heegaard Floer knot homology, was developed by Ozsv\'{a}th and Szab\'{o} and also, independently, by Jacob Rasmussen.   The original construction of the Heegaard Floer homologies required solving a system of partial differential equations.  Recently two combinatorial algorithms were developed. The first algorithm, discovered by Ciprian Manolescu, Peter Ozsv\'{a}th and Sucharit Sarkar, computes the Heegaard Floer knot homologies.

The second algorithm, discovered by Sucharit Sarkar and Jiajun Wang, computes $\HF$ and several of its knot versions including $\HFK$.  

In this paper, we discuss the Heegaard Floer knot homology $\HFK$ and bound the computational complexity of the Sarkar-Wang algorithm.  We also analyze the computational complexity of the combinatorial algorithm that calculates the Heegaard Floer homology from a nice diagram.  We start by defining the Heegaard Floer homology in section \ref{preliminaries}.  Then we analyze the Sarkar-Wang algorithm in sections \ref{SWA} and the algorithm to compute $\HFK$ from a nice diagram in \ref{sec_fromnice}.  In section \ref{modification}, we present our modifications to the Sarkar-Wang algorithm, which we believe significantly lowers the computational complexity of the algorithm.  We then examine the computational complexity of the the modified algorithm together together with the algorithm to compute $\HFK$ from a nice diagram (i.e. the process of computing $\HFK$ from a Heegaard diagram that is not necessarily nice) in section \ref{SWA_mod}. Finally, we conclude with some conjectures which, if true, would imply that the Sarkar-Wang algorithm is faster than the Manolescu-Ozsv\'{a}th-Sarkar algorithm.

We would like to thank Professor Robert Lipshitz and Thomas Peters for supervising us during the summer of 2007 in doing this research.  Many of the definitions used in this paper were adopted from handouts given to us by Robert Lipshitz.

The first author would like to thank professor Walter Neumann for supporting this research.
The second author would like to thank the I.I.Rabi Scholarship for supporting this research.
The third author would like to thank the Columbia Mathematics Summer Research Program for supporting this research.

\section{Preliminaries}\label{preliminaries}
\subsection{Heegaard Diagrams}

\ \vspace{5pt}

We begin by reviewing some notions from  \cite{Introduction} and \cite{sw} as well as introducing a some new definitions and notation.
\ \vspace{5pt}

\begin{definition}
A \textit{Heegaard diagram for} $S^3$ is a triple $(\Sigma,\alpha,\beta)$ where
\begin{enumerate}
\item $\Sigma$ is a closed, orientable surface of genus $g$ and
\item $\alpha=\{\alpha_1,\alpha_2,...,\alpha_g\}$ and $\beta=\{\beta_1,\beta_2,...,\beta_g\}$ are $g$-tuples of pairwise-disjoint embedded simple closed curves such that $\Sigma\backslash\{\alpha_1\cup\alpha_2\cup...\cup\alpha_g\}$ and $\Sigma\backslash\{\beta_1\cup\beta_2\cup...\cup\beta_g\}$ are connected and that when viewing $\Sigma$ as sitting in $S^3$ , each $\alpha_i$ bounds a disk outside $\Sigma$ and each $\beta_i$ bounds a disk inside $\Sigma$.
\end{enumerate}
\end{definition}

\begin{definition}
A \textit{doubly pointed Heegaard diagram} $\calh$ for a knot $K$ in $S^3$ is a Heegaard diagram for $S^3$ with two \textit{special} points $z, w\in\Sigma\backslash(\alpha\cup\beta)$ with the property that $K$ can be obtained by the following procedure: connect $z$ and $w$ by two paths $\gamma_1, \gamma_2$ such that $\gamma_1$ is embedded in $\Sigma\backslash\alpha$ and the path $\gamma_2$ is embedded in $\Sigma\backslash\beta$; then push $\gamma_1$ slightly out of $\Sigma$ to get a curve $\eta_1$, and push $\gamma_2$ slightly inside to get a curve $\eta_2$; then make $K$ by gluing the two curves together i.e. $K=\eta_1 \cup \eta_2$.
\end{definition}

\be{definition}{A Heegaard diagram is \textit{disk} if all regions of the Heegaard diagram are disks.}

\brem 
We are always able to make a Heegaard diagram into a disk diagram.  Sarkar-Wang give us an algorithm for this in \cite{sw} and we also present another method in section \ref{kill_nondisks} of this paper.
\erem

\bdf 
Regions with two edges are called \textit{bigons} and regions with 4 edges are called \textit{rectangles}.  The \textit{special} regions $Z$ and $W$ are the regions that contain the special points $z$ and $w$.
\edf

\begin{definition}
\begin{itemize}
\item A region is $bad$ if it is not a bigon and not a rectangle.
\item The \textit{badness} of a region $R$ with $2n$ edges is $max(n-2,0)$.  We denote the badness of a region by $b(R)$.
\item The \textit{total badness of $\calh$} is the sum of the badness of all regions other than $Z$, $\sum_{R \neq Z} b(R)$.  We denote the total badness of a diagram by $b(\calh)$.
\item  A Heegaard diagram is \textit{nice} if its total badness is 0.
\end{itemize}
\end{definition}

\bdf 
\begin{itemize}
\item
The \textit{Heegaard distance} of a region $R$ is the smallest number of $\beta$-edges a path inside $\Sigma \setminus \alpha$ going from $R$ to $Z$ would have to pass through.  The distance of a Heegaard diagram $\calh$ is the maximal distance over all bad regions of $\calh$.  
\item
The \textit{total badness of $\calh$ at distance $d$} is $b_d(\calh) \triangleq \sum^m_{n=1}b \left(D_i \right)$, where $D_1,...D_m$ are all bad regions of distance $d$.
\end{itemize}
\edf

\begin{definition}
The \textit{Heegaard} path of a region $R$ is a path inside $\Sigma \setminus \alpha$ going from $R$ to $Z$  whose number of intersections with $\beta$-edges is equal to the \textit{Heegaard distance} of the region.
\end{definition}

\bdf The $Euler\: measure$ of a region $R$ is $e(R) \triangleq \chi (R)-\frac{\#corners\: of\: R}4$, where $ \chi (R)$ is the Euler characteristic of $R$. \edf

\blem For a disk Heegaard diagram $\calh$, if  $b(\calh)$ is the total badness of $\calh$, $b_z$ is the badness of $Z$, $g$ is the genus of the diagram and $B$ is the number of bigons in a diagram, then the following equality holds
 \beq b(\calh)+b_z=4(g-1)+B \label{Euler_measure_eq} \eeq \elem

\bpf The total Euler measure of genus $g$ surface is equal to the sum of Euler measures of its regions: $$2-2g=\sum_{D} e(D)=\frac 12( -b(\calh)-b(Z)+B)$$ which becomes (\ref{Euler_measure_eq}) after rearranging.  \epf

\subsection{Heegaard Moves}
There are three moves that can be made on a Heegaard diagram $\calh$ for a knot $K$ in $S^3$ to yield another Heegaard diagram $\calh'$ for $K$ \cite{holo}:\\

\begin{itemize}

\item Isotopies: moving an $\alpha$ (or $\beta$) curve around without intersecting another, including itself, $\alpha$ (or $\beta$) curve.  Any isotopy can be obtained by a sequence of finger moves. (Figure \ref{fingermove})

\begin{figure}[h]
\center{\includegraphics[scale=.25]{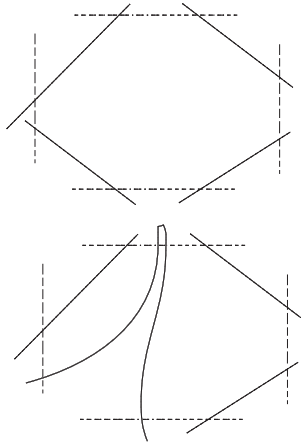}}
\caption{The dotted arcs represent $\alpha$- arcs and solid arcs represent $\beta$-arcs.  This is an example of a finger move of an $\alpha$ arc through a $\beta$ arc.}
\label{fingermove}
\end{figure}

\item Handleslides: pushing $\alpha_i$ over $\alpha_j$ or pushing $\beta_i$ over $\beta_j$. (Figure \ref{handleslide})

\begin{figure}
\center{\includegraphics[scale=.25]{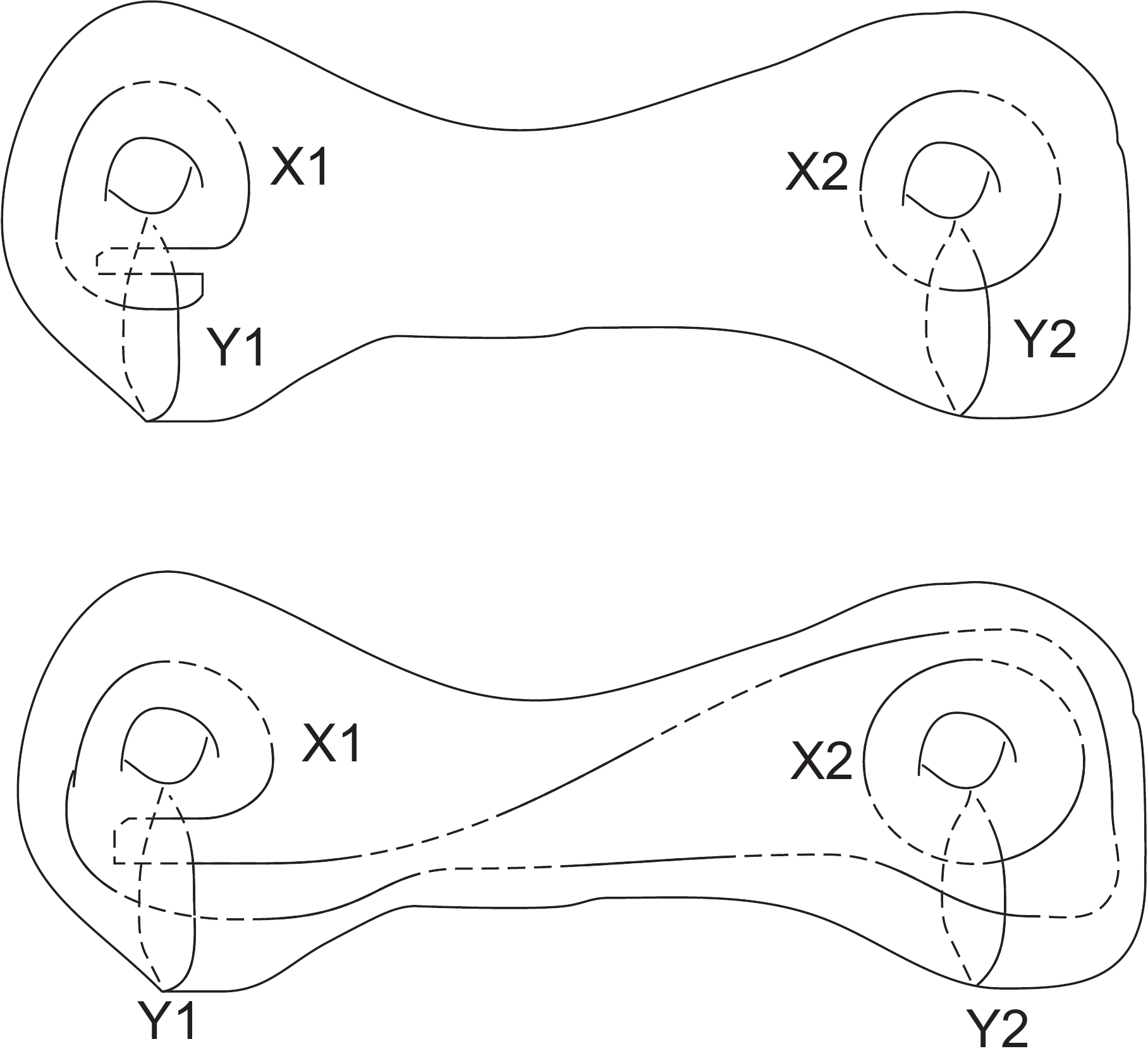}}
\caption{Handleslide: $\alpha_1$ goes around $\alpha_2$.}
\label{handleslide}
\end{figure}

\item Stabilization: gluing the Heegaard diagram to a genus 1 Heegaard diagram for $S^3$.  Destablization, which is the inverse process, is also allowed. (Figure \ref{stabilization})

\begin{figure}
\center{\includegraphics[scale=.25]{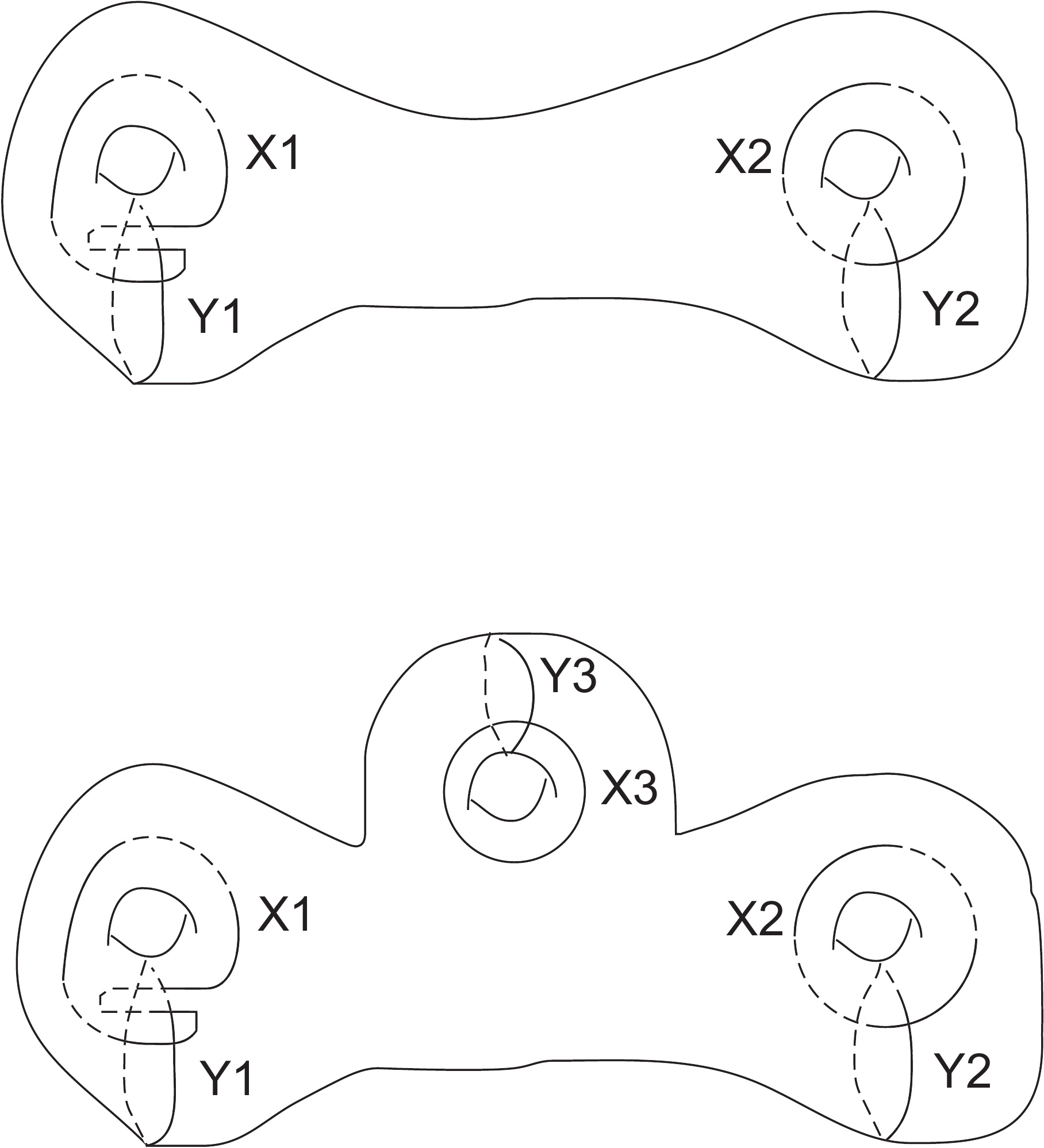}}
\caption{Stabilization of a Heegaard diagram.}
\label{stabilization}
\end{figure}

\end{itemize}

\bdf 
An \textit{alteration} of $\calh$ is any Heegaard diagram $\calh'$ obtained from $\calh$ by the the combination of the above Heegaard moves.
 \edf

\bdf
A \textit{finger move of length n} is a finger move of an $\alpha$-curve that intersects n $\beta$-curves or, respectively, a finger move of an $\beta$-curve that intersects n $\alpha$-curves.
\edf

\subsection{Calculating $\HFK$ from a Nice, Disk Heegaard Diagram}

\ \vspace{5pt}

\begin{definition}
A \textit{matching} is a $g$-tuple of points $x=\{x_1,x_2,...,x_g\} \subset(\alpha\cap\beta)$ such that exactly one $x_i$ lies on each $\alpha$ circle and exactly one on each $\beta$ circle.
\end{definition}

  In the following definition and the rest of the paper we use an unconventional, but convenient, notation.  We write $\partial_{\alpha}$ for the intersection of the image of the boundary operator on a region with the $\alpha$-curves.

\begin{definition}
A $domain$ connecting $n$-tuples $ {\bf x}=\{ x_1,...x_n \}$ and ${\bf y}=\{ y_1,...y_n \}$ is a chain of regions $D=\sum_i a_i R_i$ such that: 
	\be{itemize}{
		\item $\partial(\partial_\alpha D)=y-x$.
		\item The coefficient of $Z$ in $D$ is zero.
	}
\end{definition}

For any vertex $q$,  let $n_q(D)$ be the average of the coefficients $a_i$ of the four regions with a corner at $q$.  For a matching $\bf x$ define $n_{\bf x}(D)=\sum_{x_i \in \ \bf x} n_{x_i}(D)$.

  A graded chain complex is a sequence of homomorphisms:
  
 \begin{equation*}
\ldots \xrightarrow{\partial_{i+1}}C_{i}\xrightarrow{\partial_{i}}C_{i-1}\xrightarrow{\partial_{i-1}} \ldots  
\end{equation*}

Where the $C_i$ are abelian groups and $\partial_{i-1} \circ \partial_i =0, \hspace{5pt} \forall i$.

  A bi-graded chain complex is a collection of graded chain complexes which are arranged with each other according to a another grading.

A Heegaard Floer knot chain complex $\CFK$ s a bi-graded chain complex.  It suffices to use relative gradings in our case.  They are defined as follows: Let $D=\sum_i a_i R_i$ be the domain connecting $\bf x$ to $\bf y$.  Maslov grading $\mu$ between $x$ and $y$ is defined by $\mu ({\bf x, y})=e(D)+n_{\bf x}(D)+n_{\bf y}(D)$; and Alexander grading is defined as $A({\bf x}, {\bf y})=a_w$ where $a_w$ is the coefficient of $W$ in $D$.

The differential is a map from $CFK_{i,j}$ to $CFK_{i,j-1}$
\\
The differential is given by: $$d(\bf {x})=\sum_ {\substack{
      {\mu (\bf{ x}, \bf{ y})=1} \\
      {A \left( {\bf x}, {\bf y} \right)=0}
}}
    n_{\bf{ x}, \bf{ y}}y$$ where $n_{\bf x, \bf y}$ is the number of positive domains going between 1-	 and 2-tuples of matchings $\bf x$ and $\bf y$.
\\

Heegaard Floer homology is a bigraded complex, denoted $HFK$, whose elements are $HFK_{i,j}=Ker(d_n)/Im(d_{n+1})$

  In order to talk about the computational complexity of an algorithm we need to have some notion of time.  We deal with that now.
\subsection{A Representation of a Heegaard Diagram and a Unit of Time.}
\ \vspace{5pt}

  In our representation of a Heegaard diagram we are only interested in preserving enough information to compute $\CFK$ and $\HFK$.

Each finger move of length 1 takes $O(1)$ time to modify this presentation since we only have to add a fixed number of simplices and change the boundary data maps in finitely many places.  The modification of $BR_{i,j}$ takes $O(1)$ since the badness of at most 2 regions changed and at most one bad region is created. Similarly we need to make $O(l)$ changes for a finger move or handleslide of length $l$.

Hence, the order of the total complexity of the Sarkar-Wang algorithm is of the same order as the total length of all finger moves and handleslides, which is of the same order as the number of new 0-simplices we have added to the presentation during these finger moves and handleslides.  Hence from now on we count only the number of new 0-simplices.

\section{Sarkar-Wang Algorithm and its Complexity}
\label{SWA}

\ \vspace{5pt}
  In this section we restate and give a bound for the computational complexity of the Sarkar-Wang.

  The Sarkar-Wang algorithm works by elminating all of the bad regions at successive distances.  In order to talk about the algorithm in detail, we restate the following definitions from \cite{sw}.

  From now on $\calh$ will denote a Heegaard diagram of distance $d$ and $D_1,...,D_n$ will be its bad regions of distance $d$ such that $b \left( D_1 \right) \geq \hdots \geq b \left( D_n \right)$.

\bdf 
The \textit{distance $d$ complexity} of $\calh$ is
$$c_d(\calh)=\left(\sum_{i=1}^m b(D_i), -b(D_1), -b(D_2),\cdots, -b(D_n)\right),$$ 
and $c_d(\calh)$ is ordered lexicographically.
\edf

We will also refer to something called a \textit{chain of rectangles}.    A \textit{chain of rectangles} is set of rectangles, where each rectangle has an $\alpha$-edge in common with an another rectangle in the set and there are only two possible $\beta$-curves that the $\beta$-edges of the rectangles can be a part of.

\subsection{The Sarkar-Wang Algorithm}

\ \vspace{5pt}

  For a detailed exposition of the algorithm the reader should consult \cite{sw}.  For convenience we summarize it here.
  
\textit{Step 1: First make each $\alpha$ circle intersect at least one $\beta$ circle and vice versa.  Then eliminate the non-disk regions with finger moves. 
Step 2: Iterate the following until there are no  bad regions.  We will refer to this as the  \textbf{Sarkar-Wang procedure}:}
	
  \textit{Given a bad region $D$ of maximal distance $d$ and minimal badness among bad regions of distance $d$, Let $D_*$ be some region of distance $d-1$ adjacent to $D$ via a $\beta$ curve $\phi$.  Let $a_1, a_2, ... a_n$ be the $\alpha$ edges on the boundary of $D$ going in counterclockwise order such that $a_1$ and $a_n$ are the $\alpha$ edges adjacent to $\phi$.   Make a finger move on $\phi$ through the chain of rectangles that start from $a_2$ and ends when a non-rectangular region or region of distance less than $d$ is reached.  If the end region $E$ coincides with $D$ and the finger comes back through an $\alpha$ edge adjacent to $a_2$ (that is $a_1$ or $a_3$), then make a handle slide.} 
	
  \textit{If the end region $E$ coincides with $D$ and the finger comes back through an $\alpha$ edge other than $a_1$ or $a_3$ (we will refer to this event as the "hard" case for the rest of the paper), then make the finger move through the $a_3$ edge.  If this fails try again through the $a_4$ edge.  Continue progressing counterclockwise around the edges until you get to one (it is shown in \cite{sw} that there is one) in which the finger move ends as in one of the cases (we will refer to these cases as the "easy" cases the rest of the paper) from the previous paragraph.}

\subsection{Sarkar-Wang Lemma.}  Let us start by stating the Sarkar-Wang lemma  (lemma 4.1 in \cite{sw}).

\blem  For a distance $d$ pointed disk Heegaard diagram $\calh$ with $v(\calh)$ vertices and $c_d(\calh)\neq(0)$, by applying the Sarkar-Wang procedure we get an alteration $\calh'$, satisfying the following three properties:
\bea d(\calh')\leq d(\calh) \label{cond1}\\
c_d(\calh')<c_d(\calh) \label{cond2}\\
b(\calh') \leq b(\calh)+1 \label{cond3}\eea
\label{SW_Lemma}
\elem

Note: The original Sarkar-Wang Lemma did not contain (\ref{cond3}).  We prove it below; for  proofs of (\ref{cond1}) and (\ref{cond3}) see \cite{sw}.

\bpf[Proof of (\ref{cond3})]
	In each case the number of bigons increases by at most one, $b_z$ does not decrease and the genus stays the same.  Combining these observations with equation (\ref{Euler_measure_eq}) we get (\ref{cond3}). \epf

  This lemma shows that the Sarkar-Wang algorithm gives a nice Heegaard diagram, but it does not yield nice computational results.  To get better results we look at the algorithm more closely.
  For now we skip a detailed discussion of how step 1, the elimination of non-disk regions, affects the computational complexity and focus on the computational complexity of step 2, the Sarkar-Wang procedure.  We will look at step 1 in section \ref{kill_nondisks} where we bound the computational complexity of the entire process of calculating $\HFK$ from a Heegaard diagram.   

\blem To decrease the number of bad regions of distance $d$ we need to apply the Sarkar-Wang Procedure at most $\left(\frac{(b(D)+1)}2\right)^2$ times. The number of vertices of the diagram  increases by at most a factor of $3^{b(D)}$ \label{Kill_smallest_bad} during these applications.\elem

\bpf Let $D$ be as in the Sarkar-Wang procedure.  Let a 'reduction' stand for the process to whereby we eliminate one unit of badness from $D$ without creating any additional bad regions of distance d.

  Notice that if we consider only the ``easy" cases it would take us only $b(D)$ applications of the Sarkar-Wang procedure to eliminate the badness of the region $D$ since every time the procedure is applied, we would ``push" at least one unit of badness of $D$ into a different region.  

In the "hard" case things are a little more complicated because when we make a finger move through $a_k$ we are splitting the bad region into two bad regions $D_1$ and $D_2$.  However, this actually makes the situation easier to deal with.  For note that all of finger moves one can make in $D_2$ end up in $D_1$ (if a finger move from $D_2$ were to end up in $D_2$ then by the claim in the middle of page twelve in \cite{sw} the finger move for the "hard" case would have been made through one of the boundaries of $D_2$, a contradiction).  Therefore a reduction in the "hard" case takes at most $\left\lfloor \frac{b(D)+1}2 \right\rfloor$ applications of the Sarkar-Wang procedure.  

 Notice that during an application of a finger move we pass through each chain of rectangles at most once and so the number of vertices of the diagram have increased a factor of at most $3$.
Now $b(D_1) \leq b(D)-1$ so in order push all of the badness out of D and eliminate all of the extra bad regions of distance d, takes at most $b(D)$ 'reductions'. Hence in order to eliminate a single bad region of distance d we have to apply the Sarkar-Wang Procedure at most $\sum_{i=1}^{b(D)} \left\lfloor \frac{i+1}2 \right\rfloor \leq \left(\frac{\left(b(D)+1\right)}2\right)^2$ times and the number of vertices in the diagram increases by at most a factor of $3^{b(D)}$.\epf

\blem To eliminate all of the bad regions at a distance $d$ we need to apply Sarkar-Wang Procedure at most $\frac{(b_d+1)^2}2$ times this increases the the number of the vertices in the diagram by at most a factor of $3^{b_d(\log(b_d)+1)}$. \label{Kill_bad_at d} \elem

\bpf Since each application of lemma \ref{Kill_smallest_bad} eliminates a bad region, it takes $m$ applications of \ref{Kill_smallest_bad}, where $m$ is the number of bad regions of distance $d$.  

Note that while an application of \ref{Kill_smallest_bad} may increase the badness of a region at distance d, it does so only by reducing the badness of $D$ in the same amount.  Thus while $b_d$ might be decreased, it will never be increased by an application of \ref{Kill_smallest_bad}.  Combining this with the fact that the badness of the least bad distance $d$ bad region is at most $\frac{b_d}{i}$, we get that the the removal of a distance $d$ bad region, when there are $i$ left, takes at most $\left\lfloor \frac{(\lfloor \frac{b_d}{i} \rfloor+1)^2}4 \right\rfloor$  applications of the Sarkar-Wang procedure, and that the number of intersections of the diagram increases by a factor of at most $3^{\left\lfloor \frac{b_d}{i} \right\rfloor}$.  

Hence, by the time we have eliminated all of the bad regions of distance d, we have applied the Sarkar-Wang Procedure at most $$\sum_{i=1}^{b_d} \left\lfloor \left(\frac{\lfloor \frac{b_d}{i} \rfloor+1}2\right)^2 \right\rfloor \leq \frac{(b_d+1)^2}2$$ times and the number of intersections increased by a factor of at most $$3^{\sum_{i=1}^{b_d} \left\lfloor \frac{b_d}{i} \right\rfloor} \leq 3^{b_d (\log(b_d)+1)}.$$ \epf

With $f(x) \triangleq \frac{(x+1)^2}2$.  We iterate Lemma \ref{Kill_bad_at d} to get the following result:

\bres We can obtain a nice Heegaard diagram by applying the Sarkar-Wang Procedure at most $f(b_1+f(b_2+f(b_3+...+f(b_d)))) \leq f^d(b) \leq b^{2^d}$ times, where $b$ is the total badness and $d$ is the distance of the original Heegaard diagram.  During these applications the number of new vertices increases by at most a factor of $3^{b^{2^d}}$. \eres

  In order to improve the computational complexity of the algorithm, we introduce two modifications. 

\section{Modification of the Sarkar-Wang Procedure} \label{modification}

  We were able to shave some time off of the computational complexity of the algorithm in two areas.  The first one modifies what we have refered to as the "hard" case.  The second one has takes place at the beginning of the algorithm and makes the Heegaard diagram more amenable to the algorithm.

\begin{theorem}
We do not have to undo the finger move when the end region is $E$ and the finger move comes back through an edge other than $a_1$ or $a_3$.
\end{theorem}

\begin{proof}
In the original algorithm we continue to make and undo finger moves in a counter clockwise fashion until we find on that reduces to one of the easy cases.  But since the finger move is made through an edge $a_i$ with $i>2$  the region on the right hand side of the finger move will still be bad.  Thus we will have to redo all of the finger moves that we had previously withdrawn.  (the unmodified algorithm tells us to leave them there this time since they no longer end in the same region that they started, i.e. they end in the region on the left hand side of the finger move).

So leaving the finger moves that we were supposed to withdraw in the original algorithm does not change the final outcome.  It only saves us the time of withdrawing them. 
\end{proof}

The second modification is a more substantial change to the algorithm and also gives us a more significant reduction in computational complexity.

Recall the bound that we had computed for the number of new vertices an application of the Sarkar-Wang Algorithm resulted in: $3^{b^{2^d}}$.  Clearly this would be much improved if $d=1$.  In fact, if $d=1$ we can just apply lemma \ref{Kill_bad_at d} directly and bound the number of new vertices by $3^{b(\log(b)+1)}$.

This is the motivation for our next modification of the Sarkar-Wang algorithm: making the Heegaard diagram into a Heegaard diagram of distance 1.  This modification is more involved than the previous one so we give it its own subsection.

\subsection{Making a Heegaard Diagram into a Distance 1 Heegaard Diagram.}
\label{Dist1}

	The idea behind the modification is reduce the Heegaard distance of a region by attaching a handle to the region.  This clearly reduces the region's Heegaard distance, but in general the new $\alpha$ and $\beta$ curves may affect the Heegaard distances of other regions in the diagram.  In order to reduce the Heegaard distance of all, and not just one, of the regions of a diagram we would like the attaching of the handle to not increase the distance of any other regions.  This is the content of the next theorem.
		
	In order to keep track of the regions between $\calh$ and $\calh^{\prime}$ ($\calh$ with a handle attatched to it) we introduce a function $\iota$:  

\begin{definition}
Let $\square$ (resp. $\square^{\prime}$) be the set of regions of $\calh$ (resp. $\calh^{\prime}$).  Then define $\iota:\square^{\prime} \to \square$ to be the identity on regions that do not have part of the new $\beta$ curve as an edge and let it map the new regions (i.e. those with an edge that is part of the new $\beta$ curve) to the regions in $\calh$ whose intersection with the new $\beta$ curve gave rise to them. 
\end{definition}

\begin{theorem}
\label{Handle_thrm}
Given a region $R_*$ of distance $d>1$ and a Heegaard diagram $\calh$ there exists an alteration $\calh'$ such that:

\begin{align}
& \iota^{-1}(R_*)=\{ R_*' \}, \  \ d\left( R_*' \right)=1 \label{eq1} \\
&d(R') \leq d(R), \quad \forall R' \in \iota^{-1}(R), \forall R \in \square     \label{eq2} \\
&\sum_{R' \in \ \iota^{-1}(R)} b(R') = b(R)  \quad \forall R \in \square \setminus {R_*,Z} \label{eq3}\\
& b(\calh^\prime)  = b(\calh) +2 \label{eq4}\\
& b_{>1}(\calh')<b_{>1}(\calh) \label{eq>1},   \quad where\quad  b_{>1}=b-b_1 \quad 
 \end{align}

\end{theorem}

\begin{proof}

Construction: Attach a handle to $\mathcal{H}$ in the following manner:  Attach one end to a disk in the interior of $Z$ and the other end to a disk in the interior of $R_*$.  Draw a new $\beta$ curve around its meridian.  Construct a new $\alpha$ curve by first drawing a curve from $Z$ to $R_*'$ that goes along the handle. Continue the curve along the side of a Heegaard path of the region in the original diagram.  Gluing the ends of these two curves together we get the new $\alpha$ curve.

\begin{figure}[h]
\center{\includegraphics[scale=.5]{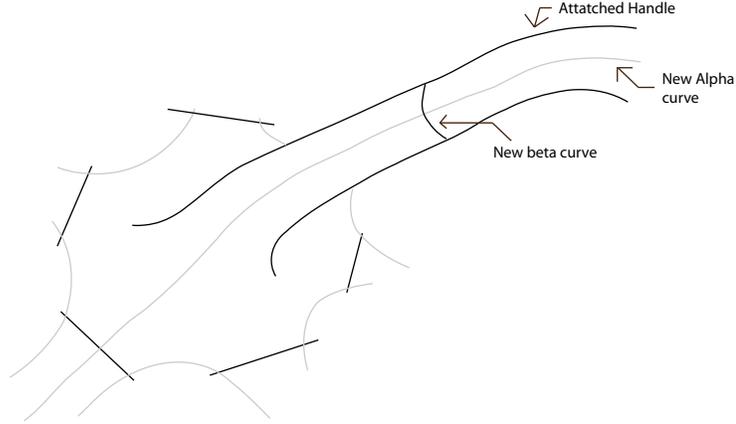}}
\caption{An example of the handle being attached to the region $R$ in our modification.}
\end{figure}

Proof of ($\ref{eq1}$): Since the new $\beta$ curve borders both $Z'$ and $R_*'$ we have ($\ref{eq1}$).

Proof of ($\ref{eq2}$): We show that, for every region in the new diagram, there is at least one Heegaard path from the original diagram that does not intersect the new $\beta$ or the new $\alpha$ curve.

\hspace{5pt}It is sufficient to consider only the regions which are split by the new $\alpha$ curve.  For if a Heegaard distance path from another region pass through one of these regions in the old diagram, from that point on it would not have gotten to $Z$ any faster than by moving along the side of a the Heegaard path of a region split by the new $\alpha$ curve.  Therefore, for any given region, there is a Heegaard path that does not intersect a Heegaard path of a region split by the new $\alpha$ curve.  Similarly, since the new $\alpha$ curve is drawn along side of a Heegaard path for every region there also exists a Heegaard path, from the old diagram, that does not intersect the new $\alpha$ curve.  Therefore the Heegaard distances of each of these regions cannot be any greater than their Heegaard distances in the original diagram.  This proves equation ($\ref{eq2}$).

Proof of ($\ref{eq3}$): Observe that if $R \in \calh^2 \setminus R_*$ is split into two regions $R_1', R_2'$ then $b(R_1')+b(R_2')=b(R)$; if $R$ is not split then $\#\{\iota^{-1}(R)\}=1$ and the regions badness does not change during the process of alteration.

Proof of ($\ref{eq4}$): Observe that the genus increased by $1$, $b(Z')=b(Z)+2$ and no new bigons were formed; combining this with equation (\ref{Euler_measure_eq}) we get the desired result.

Proof of ($\ref{eq>1}$): Combine ($\ref{eq1}$), ($\ref{eq2}$) and ($\ref{eq3}$).

\end{proof}

\bobs \label{obs_handle1}
Attaching a handle creates at most $d+1$ new vertices.
\eobs

\bcor \textit{ A Heegaard diagram of any distance can be made into a Heegaard diagram of distance 1 by stabilization.}\ecor

\begin{proof}
Apply theorem \ref{Handle_thrm} repeatedly to all of the regions.
\end{proof}

\bobs From observation \ref{obs_handle1} it follows that making a Heegaard diagram into a Heegaard diagram of distance 1 creates at most b(d+1) new vertices.
\eobs

\section{Computing the Heegaard Floer Knot Homology from a Nice Diagram}
\label{sec_fromnice}

\ \vspace{5pt}

In order to compute $\HFK$ we need to calculate a differential.  In order to compute the differential we need for any two matchings  to be able to find the domain connecting them.  We give an algorithm for finding domains here.

\subsection{Computing the ``Inverse" of $\partial\partial_{\alpha}$}

\ \vspace{5pt}  If we are able to construct the the inverse of $\partial\partial_{\alpha}$, then given any two matchings we will always be able to find the domain connecting them.  Instead it turns out to be simpler, and sufficient, to just compute the inverse of an extension of $\partial\partial_{\alpha}$ into a larger domain.

Let $\calh$ be nice Heegaard diagram of genus $g$ with $\vv$ vertices, $\ee$ edges and $\rr$ regions.  The Euler characteristics of $\Sigma$ is $2-2g=\vv-\ee+\rr$, but since the valence of each vertex is $4$ we get have $2\ee=4\vv$ and hence $$\vv=\rr+2g-2.$$
Let $X_1=\m R\langle R| R \in \calh^2, R \neq Z \rangle$ be the $\m R$-vector space generated by the regions of $\calh$ not containing $z$, i.e. the formal $\m R$-linear combinations of those regions.  Let $X_2=\m R\langle \alpha_1,\alpha_2,...\alpha_g,\beta_1,\beta_2...\beta_{g-1}\rangle$ vector space generated by $\alpha_1,\alpha_2,...\alpha_g,\beta_1,\beta_2...\beta_{g-1}$.
Let $Y=\m R(v | v \in \calh^0 )$ be vector space generated by vertices of $\calh$.

Let $f_1:X_1 \-> Y$ be the map defined by $\partial\partial_{\alpha}$ and $f_2$ be the map that assigns to each $\alpha$ or $\beta$ circle the sum of vertices that are on it (recall that $\partial_{\alpha}$ is a boundary map intersected with $\alpha$, i.e. it gives only $\alpha$ edges on the boundary of the region).

\be{claim}{The map $g \triangleq f_1 \oplus f_2:X_1 \oplus X_2 \-> Y$ is bijective.}
\bpf Since $dim(X_1)+dim(X_2)=(\rr-1)+(2g-1)=\vv=dim(Y)$, it suffices to show that the map is injective.  We show this by proving that $f_1$ and $f_2$ are injective and that $\textrm{Im}(f_1) \bigcap \textrm{Im}(f_2)=\{0\}$.

Injectivity of $f_1$ follows from theorem 5.14 in \cite{Introduction} for otherwise otherwise the domain between matchings would not be unique.

Suppose $f_2(\sum_{i=1}^n a_i \alpha_i+\sum_{i=1}^{n-1} b_i \beta_i )=0$ for some $a_i, b_i \in \m Z$.  Not including $\beta_n$ in the domain amounts to having all of the coefficients of the curves that intersect $\beta_n$ be zero (since $\beta_n$ is not there to cancel them out).  Likewise the coefficient of any curve that is connected to $\beta_n$ by some path inside $\alpha \cup \beta$ must also be zero. Since all our regions are disk regions and $\alpha \cup \beta$ is connected and thus all curves are connected to $\beta_n$.  Hence all $a_i$'s and $b_i$ are zero. Thus $f_2$ is injective.

To see $\textrm{Im}(f_1) \bigcap \textrm{Im}(f_2)=\{0\}$ first notice that we can think of $Y$ as a vector space over $\mathbb{R}$.  Give this space the inner product 
\begin{equation*}
\langle \alpha , \beta \rangle \triangleq \sum_{i=1}^{\vv}a_ib_i \quad \text{where} \quad \alpha=\sum_{v_i \in \calh^0}a_iv_i ,\beta=\sum_{v_i \in \calh^0}b_iv_i 
\end{equation*}

We now show that $\textrm{Im}(f_1) \bot \textrm{Im}(f_2)$.  It suffices to check only the basis elements of $\textrm{Im}(f_1)$ and $\textrm{Im}(f_2)$.  

Let $\alpha^0$ (resp. $\beta^0$) denote the set of points lying on a given $\alpha$ (resp. $\beta^0$) curve and $\mathfrak{C}$ be the set of all the $\alpha$ and $\beta$ circles in the diagram.  Then the basis elements of $\textrm{Im}(f_2)$ are $\{\sum_{v_i \in \alpha^0}v_i,\sum_{v_i \in \beta^0}v_i\}_{\alpha,\beta \in \mathfrak{C} }$.  

The basis elements of  $\textrm{Im}(f_1)$ are: 
\begin{equation*}
\{ v_1-v_2+v_a-v_b | v_1, v_2 (resp. v_a,v_b) \text{are on the same $\alpha$ circle and are adjacent to each other.} \}
\end{equation*}

Let $\alpha \in \textrm{Im}(f_1)$ and $\beta \in \textrm{Im}(f_2)$ then $\langle \alpha,\beta \rangle = 1-1+1-1=0$.  Thus $\textrm{Im}(f_1) \bot \textrm{Im}(f_2)$.  In particular, $\textrm{Im}(f_1) \bigcap \textrm{Im}(f_2)=\{0\}$.
\epf

Let $h=g^{-1}:Y \-> X_1 \oplus X_2$ and $h_i=\pi_i \circ h$ be its the projections on $X_i$.

Since $g$ is linear computing $h_i$ from  $g$ takes $O(v^3)$ by Gauss-Jordan Elimination.

Now if $x=(x_1,...,x_n)$ and $y=(y_1,...,y_n)$ are two $n$-tuples of vertices of $\calh$ then $h(x-y)=h((x_1+...+x_n)-(y_1+...+y_n))$ finds the domain that connects the two $n$-tuples. Hence,  since $h$ can be represented as $v \times v$ matrix, finding the domain between n-tuples takes $O(v^2)$ time.

Note that if: $h_2(x-y) \neq 0$ or $h_1(x-y)$ has non-integer entries then there is no such domain; otherwise the $h_1(x-y)$ gives a domain connecting $x$ to $y$.
\subsection{Computing the Relative Alexander and Maslov Gradings}

By examining the definitions of Alexander and Maslov gradings we see that if we know the domain connecting generator $x=(x_1,...x_g)$ to $y=(y_1,...y_g)$, then the time to compute the corresponding relative gradings is $O(v^2)$. Finding the domain connecting two generators is done by multiplication of $x-y$ by the matrix $h$, which is also $O(v^2)$.  Hence it takes $O(v^2)$ to compute the relative gradings between two generators.  Thus to compute gradings for the entire diagram, it takes $O(mv^2)$ time, where $m$ is the number of generators.

\subsection{Computing the Differential}

Now let $T$ be a free $\m Z$-module of all the 4-tuples of the vertices of $\calh$.
Construct a map $\sigma_2:T \-> \{ 0,1 \}$ which assigns to each 4-tuple $t=(t_1, t_2, t_3, t_4)$ the value 1 if the pair $((t_1,t_2),(t_3,t_4))$ bounds the positive domain not containing $z$ or $w$ and the value 0 otherwise.  Construct the corresponding map $\sigma_1$ for 2-tuples. This takes $O(v^6)$ time since there are $O(v^4)$ 4- and 2-tuples and finding the domain takes $O(v^2)$ time in each case.

Take any pair of generators ${\bf x}, {\bf y}$ of $\CFK$ such that $\mu({\bf x},{\bf y})=1$ and $A({\bf x},y)=0$. Now if $x$ and $y$ have two distinct coordinates, then ${\bf x}=(x_1,x_2,..., x_g), {\bf y}=(x_1,x_2,..., x_{g-2}, y_{g-1},  y_g)$, and $n_{{\bf x},{\bf y}}= \sigma_2 (x_{g-1}, x_g, y_{g-1}, y_g)$.  If $x$ and $y$ have one distinct coordinates, then ${\bf x}=(x_1,x_2,...x_g), {\bf y}=(x_1,x_2,...x_{g-2}, x_{g-1},  y_g)$, and $c= \sigma_2 (x_g, y_g) $.  Otherwise $n_{{\bf x},{\bf y}}=0$.

To compute each such coefficient takes $O(g^4)$ time since $\sigma_4$ and $\sigma_2$ have already been calculated and there are $O(g^4)$ 4- and 2-tuples.  Hence, computing the whole differential will take $O(m^2g^4)$ time.
\subsection{Computing the Homology and the Overall Time of the Algorithm.}
Computing the Homology is done by Gauss-Jordan Elimination and hence will take $O(m^3)$.  Thus the overall time of computing the homology from a nice Heegaard diagram is $ O(m^3+m^2g^4+v^6)$.  We can bound $m$ by $m \leq v^g$.  This gives us $O(v^{3g}+v^{2g}g^4+v^6)$.

Hence we have shown the following:

\blem Given a nice Heegaard diagram $\calh$ of genus $g$ with $v$ vertices and $m$ matchings it takes $O(m^3+m^2g^4+v^6)=O(v^{3g}+v^{2g}g^4+v^6)$ to compute its homology. \label{Homology_lemma} \elem

\section{The Modified Sarkar-Wang Algorithm}
\label{SWA_mod}

Let $\calh_0$ be a Heegaard diagram for a knot $K$.  Let $\calh_i$ be the Heegaard diagram after we complete the $i^{th}$ step of the following algorithm. Let $g_i$, $b_i$, $d_i$ and $v_i$ be respectively the genus, distance, total badness and number of vertices of $\calh_i$.  Let $t_i$ be the time that $i^{th}$ step takes.

\subsection{Step 1: Killing Bigons that do not Contain z or w.}
If we have a bigon that does not contain $z$ or $w$ and is bounded by an $\alpha$ curve $a$ and a $\beta$ curve $b$, then pull $b$ through $a$ to eliminate the bigon (think of this like a reverse finger move).  Every time this happens we decrease the number of vertices by $2$.  Hence this step is finite and it takes at most $O(\frac{v_0}2)$.

So for this step we have:

\beq g_1=g_0; \quad  v_1 \leq v_0; \quad  t_1 = O(v_0) \label{step1} \eeq

\subsection{Step 2: Killing Non-Disk Regions.}

\label{kill_nondisks}

In the algorithm found in \cite{sw} this step could possibly increase the number of bigons.  here we present a method for killing non-disk regions that leaves us with no more than two bigons.

\bdf The \textit{ugliness} of a region is the number of its boundary components minus one.  \edf

For each non-disk region $R$ make the shortest finger move that will cut $R$, i.e. go from one boundary component of $R$ to another boundary component of $R$.  Continue this finger move until it intersects one of the previous bigons.  Note that no new bigons are created in the movement except possibly at the end.  However here too, no essentially new bigons are created.  We can always intersect the old bigon since the bigon has both an $\alpha$ and a $\beta$ edge.  By intersecting the finger move with the bigon, we destroy the old one in the process of creating the new one so no essentially new bigons are created.  Each of these finger moves decreases the ugliness of $\calh$ by 1.

Now combining the following facts yield that the ugliness of $\calh_1$ is at most $2g-1$:
\be{enumerate}{
\item The Euler measure of all regions combined is $2-2g$,
\item all regions except bigons have non-positive Euler measure,
\item bigons have euler measure 1/2
\item and any non-disk region with $k+1$ boundaries has an Euler measure at most $-k$. }

Now all the finger moves can be made simultaneously and since they are shortest distance paths their length is at most $v_1$.  Hence at most $2v_1(2g_1-1)$ new vertices are created.

Lemma 2.1 gives us that  $$b(\calh_1)+b_z(\calh_1)=4(g(\calh_1)-1)+B(\calh_1)$$ where $b$ is the total badness, $b_z$ is the badness of region that contains $z$, $g$ is the genus and $B$ is the number of bigons.  Hence $b_2 \leq 4(g_1-1)+2=4g_1-2$.

So for this step we have:

\beq g_2=g_1;\quad b_2 \leq 4g_1-2;\quad  v_2 \leq 4v_1g_1; \quad  t_2 = O(v_2) \label{step1} \eeq

\subsection{Step 3: Reducing the Heegaard Distance to 1 by Attaching Handles.}
\label{reduceto1}
Attach handles as in section \ref{stabilization} until the diagram becomes distance 1.  The total number of new vertices this creates is bounded by $d_2b_2 \leq v_2b_2$.   From theorem \ref{Handle_thrm} it follows that:
\beq g_3 \leq g_2+b_2; \quad b_3 \leq 3b_2; \quad  v_3 \leq v_2b_2; \quad  t_3 = O(v_3b_2) \label{step2} \eeq

\subsection{Step 4: Apply the Modified Sarkar-Wang Procedure to Obtain a Nice Diagram.}
From result \ref{Dist1} we get
\beq g_4 \leq g_3; \quad b_4=0; \quad  v_4 \leq v_3 2^{b_3 (\log(b_3)+1)}; \quad  t_4 = O(v_4) \label{step3} \eeq

Combining these 4 steps we obtain:
\beq g_4 \leq 5g_0-2; \quad  v_4 \leq (16v_0g_0^2)2^{12g_0 (\log(12g_0)+1)}=v_0 2^{O \left( g_0\log(g_0) \right) }; \quad  \sum_{i=1}^4 t_i = O(v_4) \label{step3} \eeq

\bres Starting from a knot $K$ with a Heegaard diagram $\calh$ of genus $g$ and with $v$ vertices and applying the modified Sarkar-Wang algorithm we obtain nice Heegaard diagram $\calh'$ with at most $v'=v 2^{O(g)}$ vertices and the time to compute a nice Heegaard diagram is thus $O(v')$.  \eres

\bres Starting from a knot $K$ Heegaard diagram $\calh$ of genus $g$ and with $v$ vertices and computing the Heegaard Floer Homology through the modified Sarkar-Wang algorithm we obtain $\HFK$ in $O(v^{15g}) 2^{O(g)}$ time. \label{distance1hom_res} \eres

This comes from lemma \ref{Homology_lemma} since we know that the time to compute the homology is:
$$O(v_3^{3g_3}+v_3^{2g_3}g_3^4+v_3^6)=O((16v_0g_0^2)^{15g_0}2^{180g_0^2 (\log(12g_0)+1)})=O(v_0^{15g_0})2^{O(g_0^2 \log(g_0) )}$$.

\section{Speculations and Conclusions.}

Notice that in the result \label{distance1hom_res} for a fixed genus $g$ the algorithm is polynomial in $v$.  This leads us to some speculations below: we pose two conjectures and provide some more concrete results if they are true:

\bconj If we fix $\tau$ then for a tunnel number $\tau$ knot $K$ of crossing number $k$ we can construct a Heegaard diagram of genus $\tau+1$ for $K$ with at most $2^{O(k)}$ vertices; and the construction takes $2^{O(k)}$ time. \label{tunnel_vertex_conjecture} \econj

Combining this conjecture with result 6.2 we get:

\bres If conjecture \ref{tunnel_vertex_conjecture2} is true then for fixed $\tau$, the time it takes to compute the Heegaard Floer knot homology through the modified Sarkar-Wang algorithm is $2^{O(k)}$ where $k$ is the crossing number of the knot $K$. \eres

\bres If conjecture 1 is true, then, for fixed $\tau$, the modified Sarkar-Wang algorithm is faster than Manolescu-Ozsv\'{a}th-Sarkar algorithm for all but finitely many tunnel $\tau$ knots. \eres
\bpf Observe that Manolescu-Ozsv\'{a}th-Sarkar algorithm takes $2^{O(k \log(k))}$ to compute the Heegaard Floer Homology.  Now for fixed $t$ applying conjecture \ref{tunnel_vertex_conjecture} and result 6.2 we get that the time to compute homology is $2^{O(k)(\tau+1)} 2^{O(\tau)}$.  But $\tau$ is fixed and hence $2^{O(k)(\tau+1)}2^{O(\tau)}=2^{O(k)}$.\epf

It is also not unreasonable to assume that the following stronger conjecture is true:

\bconj If we fix $\tau$, then for a tunnel $\tau$ knot $K$ of crossing number $k$ we can construct a Heegaard diagram of genus $\tau+1$ for $K$ with at most $O(k)$ vertices; and the construction takes polynomial (in k) time. \label{tunnel_vertex_conjecture2} \econj

This would yield the following result:

\bres If conjecture \ref{tunnel_vertex_conjecture2} is true then for fixed tunnel number $\tau$, the time it takes to compute the Heegaard Floer knot homology through the modified Sarkar-Wang algorithm is polynomial in $k$, where $k$ is the crossing number of the knot $K$. \eres

If conjecture 1 is true, then the modified Sarkar-Wang algorithm is theoretically faster than the Manolescu-Ozsv\'{a}th-Sarkar algorithm.

\bibliographystyle{plain}
\bibliography{HF_references.bib}{}
\end{document}